\documentclass[12p]{article}
\usepackage{amsmath}
\usepackage{amssymb}
\usepackage{amsfonts}
\usepackage{amsthm}
\usepackage{url}
\usepackage{graphicx}

\newcommand{\E}{\ensuremath{\mathbb{E}}}

\newcommand{\Prob}{\ensuremath{\mathbb{P}}}
\newcommand{\nc}{\newcommand}
\nc{\Nset}{{\mathbb{N}}} \nc{\Rset}{{\mathbb{R}}}
\nc{\R}{{\mathbb{R}}} \nc{\N}{{\mathbb{N}}}
\nc{\Zset}{{\mathbb{Z}}}
\newcommand{\Do}{\ensuremath{\mathrm{D}[0,1]}}

\begin{document}
\newtheorem{theorem}{Theorem}[section]
\newtheorem{korollar}[theorem]{Corollary}
\newtheorem{lemma}[theorem]{Lemma}
\newtheorem{proposition}[theorem]{Proposition}

\title{\bf Analysis of radix selection on\\ Markov sources}
\author{Kevin Leckey and Ralph Neininger\thanks{Supported by DFG grant NE 828/2-1}\\
Institute for Mathematics\\
J.W.~Goethe University Frankfurt\\
60054 Frankfurt am Main\\
Germany
\and
Henning Sulzbach\thanks{Supported by the FSMP, reference: ANR-10-LABX-0098}\\
Equipe-projet RAP\\
INRIA Paris –- Rocquencourt\\
78153 le Chesnay Cedex\\
France
}

\date{April 14, 2014}
\maketitle

\begin{abstract}
 The complexity of the algorithm Radix Selection is considered for independent data generated from a Markov source. The complexity is measured by the number of bucket operations required and studied as a stochastic process indexed by the ranks; also the case of a uniformly chosen rank is considered.  The orders of mean and variance of the complexity and limit theorems are derived. We find weak convergence of the appropriately normalized complexity towards a Gaussian process with explicit mean and covariance functions (in the space $\Do$ of  c{\`a}dl{\`a}g functions on $[0,1]$ with the Skorokhod metric) for uniform data and the asymmetric Bernoulli model. For uniformly chosen ranks and uniformly distributed data the normalized complexity was known to be asymptotically normal. For a general Markov source (excluding the uniform case) we find that this complexity is less concentrated and admits a limit law with non-normal limit distribution.
\end{abstract}


\section{Introduction}
Radix Selection is an algorithm to select an order statistic from a
set of data in $[0,1]$ as follows. An integer $b\ge 2$ is fixed. In the first step the  unit interval is
decomposed into the intervals, also called {\em buckets}, $[0,1/b), [1/b,2/b), \ldots, [(b-2)/b,(b-1)/b)$ and $[(b-1)/b,1]$ and the data are assigned to these buckets according to their value. If the bucket containing the datum with rank to be selected contains further data the algorithm is recursively applied by again decomposing this bucket equidistantly and recursing. The algorithm stops once the bucket containing the rank to be selected contains no other data. Assigning a datum to a bucket is called a {\em bucket operation} and the algorithm's complexity is measured by the total number of bucket operations required.

Radix Selection is especially suitable when data are stored as  expansions in base (radix) $b$, the case $b=2$ being the most common on the level of machine data. For such expansions a bucket operation breaks down to access a digit (or bit).

In this extended abstract we study the complexity of Radix Selection in a probabilistic model. We assume that $n$ data are modeled  independently with $b$-ary expansions generated from a Markov chain on the alphabet $\{0,\ldots,b-1\}$. For the ranks to be selected we use two models. First, we consider the complexity of a random rank uniformly distributed over $\{1,\ldots, n\}$ and independent from the data. This is the model proposed and studied (for independent, uniformly over $[0,1]$ distributed data) in Mahmoud, Flajolet, Jacquet and R{\'e}gnier \cite{mafljare00}. The complexities of all ranks are averaged in this model and, in accordance with the literature, we call it the model of {\em grand averages}. Second, all possible ranks are considered simultaneously. Hence, we study the stochastic process of the complexities indexed by the ranks $1,\ldots,n$. We choose a scaling in time and space which asymptotically gives access to the complexity to select quantiles from the data, i.e., ranks of the size $t n$ with $t\in[0,1]$. We call this model for the ranks the {\em quantile-model}.

The main results of this extended abstract are on the asymptotic orders of mean and variance and limit laws for the complexity of Radix Selection for our Markov source model both for grand averages and for the quantile-model. For the quantile-model  we find Gaussian limit processes for the uniform model (defined below) for the data as well as for the asymmetric Bernoulli model (defined below). For the general Markov source model with $b=2$ we identify the first asymptotic term of the mean complexity.  For grand averages and uniform data it was shown in Mahmoud et al.~\cite{mafljare00} that  the normalized complexity is asymptotically normal.
We find that for Markov sources (with $b=2$) other than uniform  the limit distribution is no longer normal and the complexity is less concentrated.
 An explanation of this behavior is given at the end of section \ref{sec:markov:av}.

 We present our analysis separately for uniform data in section \ref{sec:uniform} and for Markov sources different from the uniform model in section \ref{sec:markov},  where the quantile-model is discussed in section \ref{sec:markov:qu}, the grand averages in section \ref{sec:markov:av}.

A general reference on bucket algorithms is Devroye \cite{de86}. A large body of probabilistic analysis of digital structures is based on methods from analytic combinatorics, see Flajolet and Sedgewick \cite{flse09}, Knuth \cite{kn98} and Szpankowski \cite{sz01}. For an approach based on renewal theory see Janson \cite{ja12} and the references given there. Our Markov source model is a special case of the model of dynamical sources, see Cl{\'e}ment, Flajolet and Vall{\'e}e \cite{clflva01}.

We close this introduction defining the Markov source model explicitly, fixing some standard notation and stating corresponding results for the related Radix Sorting algorithm.\\

\noindent
{\bf The Markov source model:} We model data strings over the alphabet $\Sigma=\{0,\ldots,b-1\}$ with a fixed integer $b\ge 2$ generated by
a homogeneous Markov chain. The data strings $s=(s_i)_{i\ge 1}$
are also interpreted as $b$-ary expansions of a real number $s\in[0,1]$ via the identification
\begin{align*}
s=\sum_{i=1}^\infty s_i b^{-i}.
\end{align*}
Conversely, if to $s\in[0,1)$ a $b$-ary expansion $s=(s_i)_{i\ge 1}$ is associated, to avoid ambiguity, we chose the expansion such that
 we have $s_i<b-1$ for infinitely many $i\in\Nset$. (For $s=1$ we use the expansion where  $s_i=b-1$ for all $i\in\Nset$.)
  The most important case is $b=2$ where the data are binary strings.

In general, a homogeneous Markov chain on $\Sigma$ is
given by its initial distribution $\mu=\sum_{\ell=0}^{b-1}\mu_\ell \delta_\ell$ on $\Sigma$
and the transition matrix $(p_{ij})_{i,j \in \Sigma}$. Here, $\delta_x$ denotes the Dirac
 measure in $x\in \R$. Hence, the initial state is $\ell$ with probability $\mu_\ell$ for $\ell=0,\ldots,b-1$. We have $\mu_\ell\in [0,1]$ and $\sum_{\ell=0}^{b-1}\mu_\ell=1$.
A transition from state $i$ to $j$ happens with probability $p_{ij}$, $i,j\in \Sigma$.
Now, a data string is generated as the sequence of states taken by the Markov chain.
In our Markov source model assumed subsequently  all data strings are independent
and  identically distributed according to the given Markov chain.

We always assume that $p_{ij}<1$ for all $i,j\in \Sigma$. Note that we do not necessarily assume the Markov chain to converge to a stationary distribution nor that it starts in a stationary distribution.

The case $p_{ij}=\mu_i=1/b$ for all $i,j\in\Sigma$  is the case where all symbols within all data are independent and uniformly distributed over $\Sigma$. Then the associated numbers are independent and uniformly distributed over $[0,1]$. We call this the {\em uniform model}. For $b=2$ the uniform model is also  called {\em symmetric Bernoulli model}. The {\em asymmetric Bernoulli model} for $b=2$ is the case where $p_{i1}=\mu_1=p$ for $i=0,1$ and a $p\in(0,1)$ with $p\neq\frac{1}{2}$. \\

\noindent
{\bf Notation.} We write $\stackrel{d}{\longrightarrow}$ for convergence in distribution and $\stackrel{d}{=}$ for equality in distribution. By $\mathrm{B}(n,p)$ with $n\in\N$ and $p\in [0,1]$ the binomial distribution is denoted, by $\mathrm{B}(p)$ the Bernoulli distribution with success probability $p$, by ${\cal N}(\mu,\sigma^2)$ the  normal distribution
with mean $\mu\in\Rset$ and variance $\sigma^2>0$. The Bachmann--Landau symbols are used.\\

\noindent
{\bf Radix Sorting.}
The Radix Sorting algorithm consists of assigning all data to the buckets as for Radix Selection. Then the algorithm recurses on all buckets containing more than one datum. Clearly, this leads to a sorting algorithm.  The complexity of Radix Sorting is also measured by the number of bucket operations. It has thoroughly been analyzed in the uniform model  with refined expansions for mean and variance involving periodic functions and a central limit law for the normalized complexity, see Knuth \cite{kn98}, Jacquet and R{\'e}gnier \cite{jare88}, Kirschenhofer, Prodinger and Szpankowski \cite{kiprsz89}  and Mahmoud et al.~\cite{mafljare00}.

For the Markov source model (with $b=2$ and $0<p_{ij}<1$ for all $i,j=1,2$) the orders of mean and variance and a central limit theorem for the complexity of Radix Sorting were derived  in Leckey, Neininger and Szpankowski \cite{lenesz13}.

\noindent
{\bf Acknowledgements.}
We thank the referees for their careful reading and constructive remarks.

\section{The uniform model --- selection of quantiles} \label{sec:uniform}
In this section our model consists of independent data, identically and uniformly distributed over $[0,1]$. We fix $b\ge 2$ and consider bucket selection using $b$ buckets in each step. The number $Y_n(\ell)$ of  bucket operations needed by bucket selection to select rank $\ell\in\{1,\ldots,n\}$ in a set of $n$ such data is studied as a process in $1\le \ell\le n$. We write $Y_n:=(Y_n(\ell))_{1\le \ell\le n}$. For a refined asymptotic analysis we normalize the process in space and time and consider  $X_n=(X_n(t))_{0\le t\le 1}$ defined for $n\ge 1$ and $t\in[0,1]$ by
\begin{align}\label{def:proc}
 X_n(t):= \frac{Y_n(\lfloor tn\rfloor+1)-\frac{b}{b-1}n}{\sqrt{n}},
\end{align}
where we set $Y_n(n+1):=Y_n(n)$. The process $X_n$ has c{\`a}dl{\`a}g paths and is considered as a random variable in $\Do$ endowed
with the Skorokhod metric $d_{sk}$, see Billingsley \cite[Chapter 3]{Billingsley1999}.

Subsequently, we use prefixes of $b$-ary expansions. For $s,t\in [0,1]$ based on their $b$-ary expansions $s=\sum_{i=1}^\infty s_i \cdot b^{-i}$, $t=\sum_{i=1}^\infty t_i \cdot b^{-i}$ with $s_i,t_i\in\{0,\ldots,b-1\}$ with the conventions stated in the introduction we denote the length of the longest common prefix by
\begin{align}\label{prefix}
j(s,t):= \max\{i\in\Nset\,|\, (s_1,\ldots,s_i)=(t_1,\ldots,t_i)\}
\end{align}
with the conventions $\max \emptyset := 0$ and $\max  \Nset := \infty$.
\begin{theorem} \label{rn_thm_in1}
Let $b\in \Nset$ with $b\ge 2$. Consider bucket selection using $b$ buckets on a set of independent data uniformly distributed on $[0,1]$. For the process $X_n=(X_n(t))_{0\le t\le 1}$ of the normalized number  of  bucket operations $Y_n(\ell)$ as defined in (\ref{def:proc}) we have weak convergence, as $n\to\infty$, in $(\Do,d_{sk})$:
 \begin{align*}
  X_n
   \stackrel{d}{\longrightarrow} G.
\end{align*}
Here, $G=(G(t))_{t\in[0,1]}$ is a centered Gaussian process (depending on $b$) with covariance function
\begin{align*}
\E[G(s)G(t)]=\frac b {(b-1)^2}- \frac{b+1}{(b-1)^2} b^{-j(s,t)},\quad s,t\in[0,1],
\end{align*}
where $j(s,t)$ is the length of the longest common prefix  defined in (\ref{prefix}) and $b^{-\infty}:=0$.
\end{theorem}
Theorem \ref{rn_thm_in1} implies the asymptotic behavior of the worst case complexity $\max_{\ell=1,\ldots,n}Y_n(\ell)$ of Radix Selection:
\begin{korollar}
For the worst case complexity of Radix Selection in the model and notation of Theorem~\ref{rn_thm_in1} we have, as $n\to\infty$, that
\begin{align*}
\frac{1}{\sqrt{n}}\left(\sup_{1\le \ell\le n}Y_n(\ell) - \frac{b}{b-1}n\right) \stackrel{d}{\longrightarrow} \sup_{t\in[0,1]} G(t).
\end{align*}
\end{korollar}

For the Gaussian process $G$ in Theorem \ref{rn_thm_in1} we have the following results on the tails of its supremum and regarding the continuity of its paths.
\begin{theorem}\label{suptails}
For the supremum $S=\sup_{t\in [0,1]}G(t)$ of the Gaussian process $G$ in Theorem \ref{rn_thm_in1} we have for any $t>0$ that
$$\Prob(|S-\E[S]|\geq t)\leq 2\exp\left(-\frac{(b-1)^2}{2b} t^2\right).$$
\end{theorem}
In the Euclidean topology on $[0,1]$ (induced by absolute value) the Gaussian process $G$ in Theorem \ref{rn_thm_in1} does not have continuous paths. Typically, in the study of Gaussian processes
a metric on the index set is derived from the covariance function. We consider
\begin{align*}
 d(s,t):=\sqrt{\E[(G(t)-G(s))^2]}=\frac{\sqrt{2(b+1)}} {b-1} \cdot b^{-j(s,t)/2},\quad s,t\in[0,1].
\end{align*}
The subsequent results in this section are stated with respect to the (topologically) equivalent metric
\begin{align*}
 d_b(s,t):=b^{-j(s,t)},\quad s,t\in[0,1].
\end{align*}
\begin{theorem}[Modulus of continuity]\label{modcont}
For the Gaussian process $G=(G(t))_{t\in[0,1]}$ in Theorem \ref{rn_thm_in1}   we have, almost surely,
$$2\frac{\sqrt{\log b}}{\sqrt{b-1}} \leq \limsup_{n\rightarrow\infty} \sup_{\substack{s,t\in[0,1],\\d_b(s,t)=b^{-n}}} \frac{ | G(t)-G(s)|}{\sqrt{n b^{-n}}}\leq  2\frac{\sqrt{2\log b}}{\sqrt{b-1}(1-b^{-1/2})}.$$
\end{theorem}
\begin{theorem}[H{\"o}lder continuity]\label{Hoelder} For any $\beta<1/2$, almost surely, the paths of the Gaussian process $G=(G(t))_{t\in[0,1]}$ in Theorem \ref{rn_thm_in1} are H{\"o}lder continuous
with exponent $\beta$ with respect to $d_b$. For any $\beta>1/2$, almost surely, the paths of $G$ are nowhere pointwise H{\"o}lder
continuous with exponent $\beta$ with respect to $d_b$.
\end{theorem}

\noindent
{\bf Outline of the analysis:}
We outline the analysis leading to Theorems \ref{rn_thm_in1}--\ref{Hoelder}.
To set up a recurrence for the process $Y_n:=(Y_n(\ell))_{1\le \ell \le n}$ we denote by $I^{n}=(I^{n}_1,\ldots,I^{n}_b)$ the numbers of elements in the $b$ buckets after distribution of all $n$ elements in the first partitioning stage. We abbreviate $F^n_0:=0$ and
\begin{align*}
F^n_r:=\sum_{j=1}^r I^{n}_j,\qquad 1\le r\le b.
\end{align*}
Note that we have $F^n_b=n$. Then, after the first partitioning phase, the element of rank $\ell$ is in bucket $r$ if and only if
$F^n_{r-1}<\ell \le F^n_{r}$. This implies the recurrence
\begin{align}\label{rec_rn1}
Y_n \stackrel{d}{=} \left( \sum_{r=1}^b {\bf 1}_{\left\{F^n_{r-1}<\ell \le F^n_{r}\right\}} Y^{r}_{I^{n}_r}\left(\ell-F^n_{r-1}\right)+n\right)_{1\le \ell\le n},
\end{align}
where $(Y^{1}_{j}),\ldots,(Y^{b}_{j}),I^{n}$ are independent and the $Y^{r}_{j}$ have the same distribution as $Y_{j}$ for all $j\ge 0$ and $r=1,\ldots,b$.

By the model of independent and uniformly distributed data we have that the vector $I^{n}$ has the multinomial $M(n;\frac{1}{b},\ldots,\frac{1}{b})$ distribution. Hence, we have
$\frac{1}{n}I^{n}\to (\frac{1}{b},\ldots,\frac{1}{b})$ almost surely as $n\to\infty$ and
\begin{align*}
\frac{I^{n}-\frac{1}{b}(n,\ldots,n)}{\sqrt{n}}\to (N_1,\ldots,N_b),
\end{align*}
where $(N_1,\ldots,N_b)$ is a multivariate normal distribution ${\cal N}(0,\Omega)$ with mean zero and covariance matrix~$\Omega$ given by
$\Omega_{ij} =\frac{b-1}{b^2}$ if $i=j$ and $\Omega_{ij} =-\frac{1}{b^2}$ if $i\neq j$.
Note that for $b=2$ we have $N_2=-N_1$.
Below, we denote by $$\mathcal{N}=(\mathcal{N}_1,\ldots,\mathcal{N}_b)$$
a vector with distribution $\frac{b}{b-1}(N_1,\ldots,N_b)$. Hence $(\mathcal{N}_1,\ldots,\mathcal{N}_b)$ has a multivariate normal distribution with mean zero and
covariance matrix $\Upsilon=(\Upsilon_{ij})_{i,j\in\Sigma}$ given by
\begin{align}\label{cov_nor_rn}
\Upsilon_{ij} = \left\{ \begin{array}{cl}
\frac{1}{b-1}, &\text{if } i=j, \vspace{3mm}\\
-\frac{1}{(b-1)^2}, &\text{if } i\neq j.
\end{array}\right.
\end{align}
For the normalized processes $X_n$ in (\ref{def:proc}) we thus obtain
\begin{align}
X_n &\stackrel{d}{=} \Bigg( \sum_{r=1}^b {\bf 1}_{\left\{F^n_{r-1}< \lfloor tn\rfloor+1 \le F^n_{r}\right\}} \sqrt{\frac{I^{n}_r}{n}}X^{r}_{I^{n}_r}
\left(\frac{nt-F^n_{r-1}}{I_r^{n}}\right)\nonumber\\
&~\qquad+
\sum_{r=1}^b {\bf 1}_{\left\{F^n_{r-1}< \lfloor tn\rfloor+1 \le F^n_{r}\right\}}\frac{b}{b-1}\frac{I^{n}_r-\frac{1}{b}n}{\sqrt{n}}
\Bigg)_{0\le t\le 1}, \label{rec_mod}
\end{align}
with conditions on independence and identical distributions as in (\ref{rec_rn1}).

To associate to recurrence  (\ref{rec_mod}) a limit equation in the spirit of the contraction method we introduce the indicator functions
\begin{align*}
{\bf I}_r(x):= {\bf 1}_{[r-1,r)}(x) \text{ for } r=1,\ldots,b-1,\qquad
{\bf I}_b(x):= {\bf 1}_{[b-1,b]}(x)
\end{align*}
and the sawtooth function $s_b:[0,1]\to [0,1]$
\begin{align*}
s_b(t):=\left\{\begin{array}{cl}
bt-\lfloor bt \rfloor, & 0\le t<1,\\
1,& t=1.\end{array} \right.
\end{align*}
Moreover, we use the  transformations $\mathfrak{A}_r:\Do\to\Do$ for $r=1,\ldots, b$ with
\begin{align*}
f\mapsto \mathfrak{A}_r(f), \quad
\mathfrak{A}_r(f)(t)= {\bf I}_r(tb)f(s_b(t)) \text{ for } t\in[0,1],
\end{align*}
and $\mathfrak{B}:\Rset^b\to\Do$ with (for $v=(v_1,\ldots,v_b)$)
\begin{align*}
v\mapsto \mathfrak{B}(v),\quad
\mathfrak{B}(v)(t)= \sum_{r=1}^b{\bf I}_r(tb)v_r \text{ for } t\in[0,1].
\end{align*}
Then we associate the limit equation
\begin{align}\label{fpe_rn}
X\stackrel{d}{=}\sum_{r=1}^b \frac{1}{\sqrt{b}}\mathfrak{A}_r(X^{r})+\mathfrak{B}(\mathcal{N}),
\end{align}
where $X^{1},\ldots,X^{b}, \mathcal{N}$ are independent, the $X^{r}$ are identically distributed  random variables with values in $(\Do,d_{sk})$ and  distribution of $X$, and $\mathcal{N}$ has the centered multivariate normal distribution with covariance matrix $\Upsilon$ given in (\ref{cov_nor_rn}).

A distributional fixed-point equation related to (\ref{fpe_rn}) appeared in Sulzbach et al.~\cite{sunedr14}, see the map $T$ in equation (2.5) of \cite{sunedr14}. The proof of our Theorem \ref{rn_thm_in1} can  be carried out analogously to the proof in sections 2.1, 2.2 and 3.3 of  \cite{sunedr14}. Note that in analogy to Lemma 2.3 in \cite{sunedr14} our fixed-point equation~(\ref{fpe_rn}) characterizes the Gaussian limit process $G$ in Theorem~\ref{rn_thm_in1} as the unique fixed-point of (\ref{fpe_rn}) subject to the constraint $\E[\|X\|_\infty^{2+\varepsilon}]<\infty$ for any $\varepsilon>0$ and, hence, in the analysis one has to adapt exponents appropriately. The proofs of our Theorems \ref{suptails}--\ref{Hoelder} can be carried out as the corresponding results in section 4 of \cite{sunedr14} which are related to and partly based on fundamental work on Gaussian processes, see Dudley \cite{dudley73}, Talagrand \cite{tala87}, Adler \cite{adler90} and
Boucheron, Lugosi and Massart \cite{boluma13}.

\section{The Markov source model}\label{sec:markov}
Now the Markov source model is considered for the data. For the rank to be selected the quantile-model is studied in section \ref{sec:markov:qu}, the model of grand averages in section \ref{sec:markov:av}. We restrict ourselves to the study of Radix Selection using $b=2$ buckets.

\subsection{Selection of quantiles}\label{sec:markov:qu}
We  consider  the complexity of Radix Selection with $b=2$ buckets assuming the Markov source model for the data and the quantile-model  for the rank to be selected.
We first define functions $m_\mu:[0,1]\to (0,\infty)$ which appear in the average complexity. For $n\ge 1$ and $i=0,1$ we recursively define sets ${\cal D}_n^i=\{s^i_{n,k}\,|\,k=0,\ldots,2^n\}$ as follows: For $n=1$ we set $(s^i_{1,0},s^i_{1,1},s^i_{1,2}):=(0,p_{i0},1)$ for $i=0,1$. Further, for all $n\ge 1$, $i=0,1$ and $0\le k\le 2^n$ we set
\begin{align*}
s^i_{n+1,k}:=\left\{\begin{array}{cl}
s^i_{n,k/2}, &\text{if } k \text{ mod } 4 \in \{0,2\},\\
p_{00}s^i_{n,(k+1)/2}+p_{01}s^i_{n,(k-1)/2},  &\text{if } k\text{ mod } 4 =1.\\
p_{10}s^i_{n,(k+1)/2}+p_{11}s^i_{n,(k-1)/2},  &\text{if } k \text{ mod } 4 =3,
\end{array}
\right.
\end{align*}
We further define ${\cal D}^i_\infty:=\cup_{n=1}^\infty {\cal D}^i_n$. Note that for each $n\ge 1$ the set ${\cal D}^i_n$ decomposes the unit interval into $2^n$ sub-intervals. For $t\in[0,1]\setminus {\cal D}^i_\infty$ we denote by
$\lambda^i_n(t)$ the length of the (unique) sub-interval of this decomposition that contains $t$. Then, for $i=0,1$ and $t\in[0,1]\setminus {\cal D}^i_\infty$  we set
\begin{align*}
m_i(t)&:=1+\sum_{n=1}^\infty \lambda_n^i(t).
\end{align*}
Further, for an initial distribution $\mu=\mu_0\delta_0 + \mu_1\delta_1$  with $\mu_0\in[0,1]$ we denote
\begin{align*}
{\cal D}^\mu_\infty:= \mu_0{\cal D}^0_\infty \cup \left(\mu_0+\mu_1{\cal D}^1_\infty\right)
\end{align*}
and, for $t\in [0,1]\setminus {\cal D}^\mu_\infty$,
\begin{align*}
m_\mu(t):=\left\{ \begin{array}{cl}
\mu_0 m_0\left(\frac{t}{\mu_0}\right)+1,&\text{if } t<\mu_0,\\
(1-\mu_0) m_1\left(\frac{t-\mu_0}{1-\mu_0}\right)+1,&\text{if } t>\mu_0.
\end{array}\right.
\end{align*}

We have the following asymptotic behavior of the average complexity:
\begin{theorem}\label{theoqmmm}
Let $Y^\mu_n(\ell)$ denote the number of bucket operations  of Radix Selection with $b=2$ selecting a  rank $1\le \ell\le n$ among  $n$ independent data  generated from the  Markov source model with initial distribution $\mu=\mu_0\delta_0 + \mu_1\delta_1$  where $\mu_0\in[0,1]$ and transition matrix $(p_{ij})_{i,j\in\{0,1\}}$ with $p_{ij}<1$ for all $i,j=0,1$. Then, for all $t\in[0,1]\setminus D^\mu_\infty$ as $n\to\infty$, we have
\begin{align}\label{rn_0904}
\E[Y^\mu_n(\lfloor tn\rfloor +1)]=m_\mu(t) n +o(n).
\end{align}
\end{theorem}

\noindent
{\bf Outline of the analysis:}
We denote by $Y^0_n=(Y^0_n(\ell))_{1\le \ell\le n}$ and $Y^1_n=(Y^1_n(\ell))_{1\le \ell\le n}$ the number of bucket operations
for a Markov source model as in Theorem~\ref{theoqmmm}  for initial distributions $p_{00}\delta_0+p_{01}\delta_1$ and $p_{10}\delta_0+p_{11}\delta_1$ respectively. Then we have the system  of recursive distributional equations, for $n\ge 2$,
\begin{align}\label{sys_yin}
Y^i_n &\stackrel{d}{=} \left({\bf 1}_{\{\ell\le J^i_n\}} Y^0_{J^i_n}(\ell) +
{\bf 1}_{\{\ell>J^i_n\}} Y^1_{n-J^i_n}(\ell-J^i_n) + n\right)_{1\le \ell\le n},\quad i=0,1,
\end{align}
where $Y^0_0,\ldots,Y^0_n,Y^1_0,\ldots,Y^1_n, J^0_n,J^1_n$ are independent
(the independence between $J^0_n$ and $J^1_n$ is not required) and we  have that $J^i_n$ is $\mathrm{B}(n,p_{i0})$ distributed for $i=0,1$.
Moreover, for general initial distribution $\mu$ we further have
\begin{align}\label{sys_yin_1}
Y^\mu_n &\stackrel{d}{=} \left({\bf 1}_{\{\ell\le K_n\}} Y^0_{K_n}(\ell) +
{\bf 1}_{\{\ell>K_n\}} Y^1_{n-K_n}(\ell-K_n) + n\right)_{1\le \ell\le n},
\end{align}
where $Y^0_0,\ldots,Y^0_n,Y^1_0,\ldots,Y^1_n, K_n$ are independent and
$K_n$ has the binomial $\mathrm{B}(n,\mu_0)$ distribution.

The proof of Theorem~\ref{theoqmmm} is based on $k$ times iterating the system~(\ref{sys_yin})  with $k=k(n)=\Theta(\log n)$ chosen appropriately. The contributions of the toll functions within these $k$ iterations yield the main contribution, the other terms are asymptotically negligible.\\

A distributional analysis of the quantile-model is left for the full paper version of this extended abstract as well as the behavior at the $t\in {\cal D}^\mu_\infty$. For these $t$ the expansion (\ref{rn_0904}) still holds when defining $m_\mu(t)$ as the average of the left-hand and right-hand limit of $m_\mu$ at $t$.

A special case where the analysis is simplified considerably is the asymmetric Bernoulli  model discussed next where we obtain a functional limit law as for the uniform model in Theorem~\ref{rn_thm_in1}. \\

\noindent
{\bf The asymmetric Bernoulli model:} The data model called  {\em asymmetric Bernoulli model} consists of all data being independent and having independent bits all identically distributed over $\Sigma=\{0,1\}$ with Bernoulli $\mathrm{B}(p)$ distribution for a fixed $p\in(0,1)$ with $p\neq \frac{1}{2}$. Note that this can also be considered as a special case of the Markov source model by choosing $\mu_0=p_{00}=p_{10}=1-p$ and $\mu_1=p_{01}=p_{11}=p$. Here, the analysis simplifies considerably compared to the general Markov source model due to the fact that the mean function $m$ corresponding to $m_0,m_1,m_\mu$ in Theorem~ \ref{theoqmmm} becomes an affine function. We have the following results:
\begin{theorem}\label{thm0201}
Consider bucket selection using $b=2$ buckets on a set of $n$ independent data  generated from the  asymmetric Bernoulli model with success probability $p\in(0,1)\setminus\{\frac{1}{2}\}$. For the process $X^\mathrm{asyB}_n=(X^\mathrm{asyB}_n(t))_{0\le t\le 1}$ of the normalized number  of  bucket operations $Y^\mathrm{asyB}_n(\ell)$ defined by
\begin{align*}
X^\mathrm{asyB}_n(t):=\frac{Y^\mathrm{asyB}_n(\lfloor tn\rfloor +1)-m(t)n}{\sqrt{n}},\quad t\in[0,1],
\end{align*}
with $Y^\mathrm{asyB}_n(n+1):=Y^\mathrm{asyB}_n(n)$ and
\begin{align*}
m(t)=\frac{2p-1}{p(1-p)}t+\frac{1}{p},\quad t\in[0,1],
\end{align*}
we have weak convergence, as $n\to\infty$, in $(\Do,d_{sk})$:
 \begin{align*}
  X^\mathrm{asyB}_n
   \stackrel{d}{\longrightarrow} G^\mathrm{asyB}.
\end{align*}
Here, $G^\mathrm{asyB}=(G^\mathrm{asyB}(t))_{t\in[0,1]}$ is a centered Gaussian process (depending on $p$) with covariance function given, for $s,t\in[0,1]$,  by
\begin{align*}
\E[G^\mathrm{asyB}(s)G^\mathrm{asyB}(t)]=-\prod_{k=1}^{r(s,t)} p[g(t,k)] + \sum_{k=1}^{r(s,t)} \frac{\prod_{j=1}^{k} p[g(t,j)]}{p[1-g(t,k)]},
\end{align*}
where $p[0]:=1-p$, $p[1]:=p$ and the functions $r:[0,1]^2\to \Nset_0\cup\{\infty\}$ and $g: [0,1] \times \N_0 \rightarrow \{0,1\}$ are defined as follows:
\begin{align*}
r(s,t)=\max\{ n\in\N_0 | g(s,\ell)=g(t,\ell),\; 1\leq \ell \leq n\}
\end{align*}
and $g$ and $h: [0,1]\times \N_0 \rightarrow [0,1]$ are recursively defined by  $g(t, 0)=0$, $h(t,0)=t$ for $t\in[0,1]$ and for $k\ge 1$ by
\begin{align*}g(t,k)&=\left\{\begin{array}{cl} 0,&
\text{if } h(t,k-1) < 1-p,\\
1,& \text{if } h(t,k-1)\geq 1-p,
         \end{array}\right.\\
 h(t,k)&=\left\{\begin{array}{cl} \frac{h(t,k-1)}{1-p},& \text{if } h(t,k-1) < 1-p,\\
          \frac{h(t,k-1)-(1-p)}{p},& \text{if } h(t,k-1)\geq 1-p.
         \end{array}\right.
\end{align*}
\end{theorem}

For the maximum of the complexities we obtain the following corollary:
\begin{korollar}
In the model and notation of Theorem~\ref{thm0201} we have, as $n\to\infty$, that
\begin{align*}
\frac{1}{\sqrt{n}}\left(\sup_{1\le \ell\le n}\left(Y^\mathrm{asyB}_n(\ell) - m\left(\frac{\ell}{n}\right)n\right)\right) \stackrel{d}{\longrightarrow} \sup_{t\in[0,1]} G^\mathrm{asyB}(t).
\end{align*}
\end{korollar}
\begin{theorem}
For the supremum $S'=\sup_{t\in [0,1]}G^\mathrm{asyB}(t)$ of the Gaussian process $G^\mathrm{asyB}$ in Theorem~\ref{thm0201} we have for any $t>0$ with $p_\vee:=\max\{p,1-p\}$ that
$$\Prob(|S'-\E[S']|\geq t)\leq 2\exp\left(-\frac{(1-p_\vee)^2}{2p_\vee} t^2\right).$$
\end{theorem}

\subsection{Selection of a uniform rank}\label{sec:markov:av}
We now consider  the complexity of Radix Selection with $b=2$ buckets assuming the Markov source model for the data and the model of grand averages for the rank.
We have the following asymptotic behavior:
\begin{theorem}\label{theorrmm}
Let $W_n$ denote the number of bucket operations  of Radix Selection with $b=2$ selecting a uniformly distributed rank independent from  $n$ independent data  generated from the  Markov source model with initial distribution $\mu=\mu_0\delta_0 + \mu_1\delta_1$  where $\mu_0\in[0,1]$ and transition matrix $(p_{ij})_{i,j\in\{0,1\}}$ with $p_{ij}<1$ for all $i,j=0,1$. Then, as $n\to\infty$,
we have
\begin{align*}
\E[W_n]=\kappa_\mu n +o(n)
\end{align*}
with $\kappa_\mu >0$
given in (\ref{eq_kappa})
and
\begin{align*}
\frac{W_n}{n} \stackrel{d}{\longrightarrow} Z_\mu,
\end{align*}
where the convergence also holds with all moments. The distribution of $Z_\mu$ is given by
\begin{align}\label{nn3101}
Z_\mu\stackrel{d}{=} B_{\mu_0}\mu_0Z^0 + (1-B_{\mu_0})(1-\mu_0)Z^1+1,
\end{align}
where $B_{\mu_0},Z^0,Z^1$ are independent and $B_{\mu_0}$ has the Bernoulli distribution $\mathrm{B}(\mu_0)$. The distributions of $Z^0$ and $Z^1$ are the unique integrable solutions of the system (\ref{lim_gams}).
\end{theorem}
\noindent
{\bf Outline of the analysis:} We denote by $W_n^\mu:=W_n$ the complexity as stated in Theorem~\ref{theorrmm} and, for initial distributions $p_{i0}\delta_0+p_{i1}\delta_1$, write $W_n^i:=W_n^{p_{i0}\delta_0+p_{i1}\delta_1}$ for $i=0,1$. We have the system of distributional recurrences, for $n\ge 2$,
\begin{align}\label{rec_gams}
W^i_n &\stackrel{d}{=} B_{ii}W^0_{J^i_n} + (1-B_{ii})W^1_{n-J^i_n} +n, \quad i=0,1,
\end{align}
where $W^0_1,\ldots,W^0_n,W^1_1,\ldots,W^1_n$ and $(J^0_n,J^1_n,B_{00},B_{11})$ are independent and we have that $J^i_n$ is Binomial $\mathrm{B}(n,p_{i0})$ distributed and $B_{ii}$ is mixed Bernoulli distributed with distribution $\mathrm{B}(J^i_n/n)$ for $i=0,1$. We normalize
\begin{align*}
Z_n^i:=\frac{W^i}{n},\quad n\ge 1,\;i=0,1
\end{align*}
and obtain, for all $n\ge 2$ that
\begin{align*}
Z^i_n &\stackrel{d}{=} B_{ii}\frac{J^i_n}{n}Z^0_{J^i_n} + (1-B_{ii})\frac{n-J^i_n}{n}Z^1_{n-J^i_n} +1, \quad i=0,1,
\end{align*}
with independence relations as in (\ref{rec_gams}). This leads to the limit system
\begin{align}\label{lim_gams}
Z^i &\stackrel{d}{=} B_{p_{i0}}p_{i0}Z^0 + (1-B_{p_{i0}})(1-p_{i0})Z^1 +1, \quad i=0,1,
\end{align}
where $Z^0,Z^1$ and $B_{p_{i0}}$ are independent and $B_{p_{i0}}$ has the Bernoulli $\mathrm{B}(p_{i0})$ distribution for $i=0,1$. It is easy to show that subject to $\E[|Z^i|]<\infty$ the limit system (\ref{lim_gams}) has a unique solution; cf.~Knape and Neininger \cite[section 5]{knne13}. Also, the convergences $Z^i_n\to Z^i$ can be shown by a contraction argument in any Wasserstein $\ell_p$ metric with $p\ge 1$. From the limit system (\ref{lim_gams}) we obtain for the expectations $\kappa_i:=\E[Z^i]$ for $i=0,1$ that
\begin{align}
\kappa_0 &= \frac{1+p_{01}^2-p_{11}^2}{2(p_{00}+p_{11})(1+p_{00}p_{11})-2(p_{00}+p_{11})^2}>0,
\label{def_ka0}\\
\kappa_1 &=
\frac{1+p_{10}^2-p_{00}^2}{2(p_{00}+p_{11})(1+p_{00}p_{11})-2(p_{00}+p_{11})^2}>0
.\label{def_ka1}
\end{align}

Now, for a general initial distribution $\mu$ we have
\begin{align*}
W_n^\mu\stackrel{d}{=} B_{\mu\mu}W_{K_n}^0 + (1-B_{\mu\mu})W_{n-K_n}^1+n,
\end{align*}
where $W^0_1,\ldots,W^0_n,W^1_1,\ldots,W^1_n,(K_n,B_{\mu\mu})$ are independent, $K_n$ has the binomial $\mathrm{B}(n,\mu_0)$ distribution and $B_{\mu\mu}$ has the mixed Bernoulli $\mathrm{B}(K_n/n)$ distribution. This implies for the limit $Z_\mu$ of $W^\mu_n/n$ the  representation
\begin{align*}
Z_\mu\stackrel{d}{=} B_{\mu_0}\mu_0Z^0 + (1-B_{\mu_0})(1-\mu_0)Z^1+1,
\end{align*}
where $B_{\mu_0},Z^0,Z^1$ are independent and $B_{\mu_0}$ has the Bernoulli distribution $\mathrm{B}(\mu_0)$. Hence, we obtain for $\kappa_\mu:=\E[Z_\mu]$ the representation
\begin{align}\label{eq_kappa}
\kappa_\mu= \mu_0^2\kappa_0+(1-\mu_0)^2\kappa_1+1,
\end{align}
with $\kappa_0,\kappa_1$ given in (\ref{def_ka0}) and (\ref{def_ka1}).
The claims of Theorem~\ref{theorrmm} follow from this outline by an application of the contraction method within the Wasserstein metrics. \\

\noindent
{\bf A remark on concentration for grand averages:} Note that for the special case $p_{ij}=\mu_i=\frac{1}{2}$ for $i=0,1$ the Markov source model reduces to the uniform model. In this case it was shown (together with more refined results) in Mahmoud et al.~\cite{mafljare00} that the complexity $W_n$ in Theorem~\ref{theorrmm} as $n\to\infty$ satisfies
\begin{align}\label{nn31012}
\frac{W_n-2n}{\sqrt{2n}} \stackrel{d}{\longrightarrow} {\cal N}(0,1).
\end{align}
Our Theorem~\ref{theorrmm} also applies: We find that for the uniform model  system (\ref{lim_gams}) is solved deterministically by $Z^i=2$ almost surely for $i=0,1$, hence plugging in into (\ref{nn3101}) we also obtain $Z_\mu=2$ and thus $W_n/n\to 2$, which is, although  only a law of large numbers, consistent with (\ref{nn31012}). (The full limit law in (\ref{nn31012}) is a corollary to our Theorem \ref{rn_thm_in1}.)

However, for $(p_{00},p_{01},p_{10},p_{11})\neq(\frac{1}{2},\frac{1}{2},\frac{1}{2},\frac{1}{2})$
 the system (\ref{lim_gams}) does no longer solve deterministically, so that $W_n/n$ then has a nondeterministic limit and is less concentrated, typical fluctuations are of linear order compared to $\sqrt{n}$ for the uniform model. This behavior becomes transparent when looking at the quantile-model: In the uniform model we have the same leading linear term $2n$ in the expansion of the means of all quantiles, which  implies that a uniformly chosen rank conditional on its size will always lead to a complexity of the same linear order $2n$. This does no longer hold for a non-uniform Markov model. The constant $m_\mu(t)$ in the linear growth $m_\mu(t)n$ of the complexity depends on the quantiles $t\in[0,1]$, see Theorem~\ref{theoqmmm}. This implies that the complexity can no longer remain concentrated: The fluctuations are now forced to be at least of linear order since different choices of the ranks lead to different linear orders. This is consistent with the fact that in Theorem~\ref{theorrmm} we then find a non-deterministic limit for $W_n/n$ and a variance of the order $\Theta(n^2)$. Further note, that this also implies a simple representation of the limit distribution of $Z_\mu$ in Theorem~\ref{theorrmm} as
 \begin{align*}
 Z_\mu \stackrel{d}{=} m_\mu(U),
 \end{align*}
 with $m_\mu$ as in Theorem~\ref{theoqmmm} und $U$ uniformly distributed on $[0,1]$.

\end{document}